\documentclass[conference]{IEEEtran}
\usepackage{fullpage,graphicx,graphics,psfrag,amsmath,amsfonts,verbatim}
\usepackage[small,bf]{caption}
\usepackage{url,ar}
\usepackage{amssymb,enumitem}
\usepackage{tikz}
\usepackage{hyperref}
\usetikzlibrary{positioning,calc}
\definecolor {processblue}{cmyk}{0.96,0,0,0}

\usepackage{fullpage,url,verbatim}

\newcommand{\reals}{{\mbox{\bf R}}}
\newcommand{\eg}{{\it e.g.}}

\newcommand{\etc}{{\it etc.}}

\title{A New Architecture for\\Optimization Modeling Frameworks}

\author{\IEEEauthorblockN{Matt Wytock, Steven Diamond, Felix Heide, Stephen Boyd}
\IEEEauthorblockA{Department of Electrical Engineering\\
Stanford University\\
\{mwytock, diamond, fheide, boyd\}@stanford.edu}
}

\IEEEoverridecommandlockouts

\IEEEpubid{\begin{minipage}{\textwidth}\ \\[12pt]
  PyHPC2016; Salt Lake City, Utah, USA; November 2016\\
  978-1-5090-5220-2/16/\$31.00 \copyright 2016 IEEE
\end{minipage}}

\begin{document}


\maketitle

\begin{abstract}
  We propose a new architecture for optimization modeling frameworks
  in which solvers are expressed as computation graphs in a framework
  like TensorFlow rather than as standalone programs built on a
  low-level linear algebra interface.
  Our new architecture makes it easy for modeling frameworks to support
  high performance computational platforms like GPUs and distributed clusters,
  as well as to generate solvers specialized to individual problems.
  Our approach is particularly well adapted to first-order and indirect
  optimization algorithms.
  We introduce \texttt{cvxflow}, an open-source convex optimization
  modeling framework in Python based on the ideas in this paper,
  and show that it outperforms the state of the art.
\end{abstract}

\IEEEpeerreviewmaketitle

\section{Introduction}

Optimization offers a principled approach to solving problems in a wide
variety of application domains, such as machine learning, statistics, control,
signal and image processing, networking, engineering design, finance, and many
others \cite{BoV:04}. Instead of
designing specialized algorithms for each individual problem, the user
describes the problem as the minimization of a cost function and the optimal
solution with minimal cost is found by the optimization method.

The wealth of
applications for this methodology has driven the development of several
high-level modeling languages. These languages provide a separation of concerns
between the development of mathematical models and the implementation of numerical
routines to optimize these models. This is especially useful for rapidly prototyping
new applications, allowing practitioners to easily experiment with different
cost functions and constraints by writing expressions that closely mimic the
mathematical optimization problem. Prominent examples of modeling languages
and frameworks include AMPL \cite{AMPL}, YALMIP \cite{YALMIP}, and CVX \cite{cvx},
as well as several tied closely to particular solvers, such as CPLEX's ILOG
\cite{cplex200711} and MathProg from GLPK \cite{makhorin2000modeling}.

Despite the popularity of these modeling frameworks,
support for modern large-scale computational environments such as GPUs and
distributed clusters is virtually nonexistent. In part, this is due to
fundamental challenges in scaling interior point methods, which have
historically been the basis for solvers of most modeling frameworks, as these methods require solving
sparse linear systems to high accuracy and as such do not benefit greatly
from GPU implementation. In addition, distributing such methods beyond a single
machine typically requires high bandwidth interconnects such as those available
exclusively in HPC environments.

However, there are also highly practical reasons for the lack of
support for new environments: mature solvers often require several years to
develop and writing entirely new implementations of low-level numerical
routines specialized to each environment is unappealing. Traditionally, a degree
of platform independence has been provided by implementing on top of low-level
linear algebra libraries (\eg, BLAS, LAPACK, and SuiteSparse),
but as we discuss in this paper, this architecture is often insufficient,
especially for large problems. In addition, such libraries do not handle memory
management and data transfer between GPU and CPU or between multiple machines.

The solution that we explore is a new architecture for
optimization modeling frameworks based on solvers represented as computation
graphs. This architecture is well-suited for solving
large optimization problems by taking advantage of problem-specific structure. In
particular, the computation graph abstraction naturally represents the
composition of structured linear operators which can be significantly
more efficient than the standard sparse or dense matrix representation. We develop such
a method in this paper and demonstrate that it outperforms the
existing state of the art for solving large convex optimization problems, a
GPU-enabled version of SCS \cite{SCSpaper}, which itself is one of the few
GPU-optimized solvers available, POGS \cite{fougner2015parameter} being
another example.

A secondary, but not insignificant, benefit of this approach is automatic support
for a wide variety of computational environments (CPU, GPU, distributed
clusters, \etc), leveraging the considerable momentum and engineering effort of
existing computation graph frameworks from the deep learning community. A
potential drawback of our approach is that the
runtime system must support the necessary mathematical operations to implement
numerical optimization algorithms. For first-order and indirect solvers, the
many frameworks developed for deep learning, such as TensorFlow
\cite{tensorflow2015}, Theano \cite{B+etal:2010-scipy,Bastien-Theano-2012},
Caffe \cite{jia2014caffe}, and Torch \cite{Torch}, provide all the necessary
functionality. The frameworks have only limited support, however, for the sparse
matrix factorization routines used by direct solvers. Thus, given the
computation graph implementations available at this time, our architecture tends to favor first-order
and indirect methods as opposed to interior point methods.



The outline of the paper is as follows.
In \S\ref{sec-traditional}, we review the traditional architecture for
optimization modeling frameworks and discuss its shortcomings.
In \S\ref{sec-alternative}, we explore prior work that addressed
the shortcomings of the traditional architecture.
In \S\ref{sec-proposed}, we describe the new architecture we propose
and the computation graph abstraction the architecture is based on.
In \S\ref{sec-examples}, we present \texttt{cvxflow}, an open-source implementation
of the ideas in this paper, and numerical results comparing \texttt{cvxflow}
with the state of the art.

\section{Traditional architecture}\label{sec-traditional}

\begin{figure}
  \centering
  \includegraphics[scale=0.5]{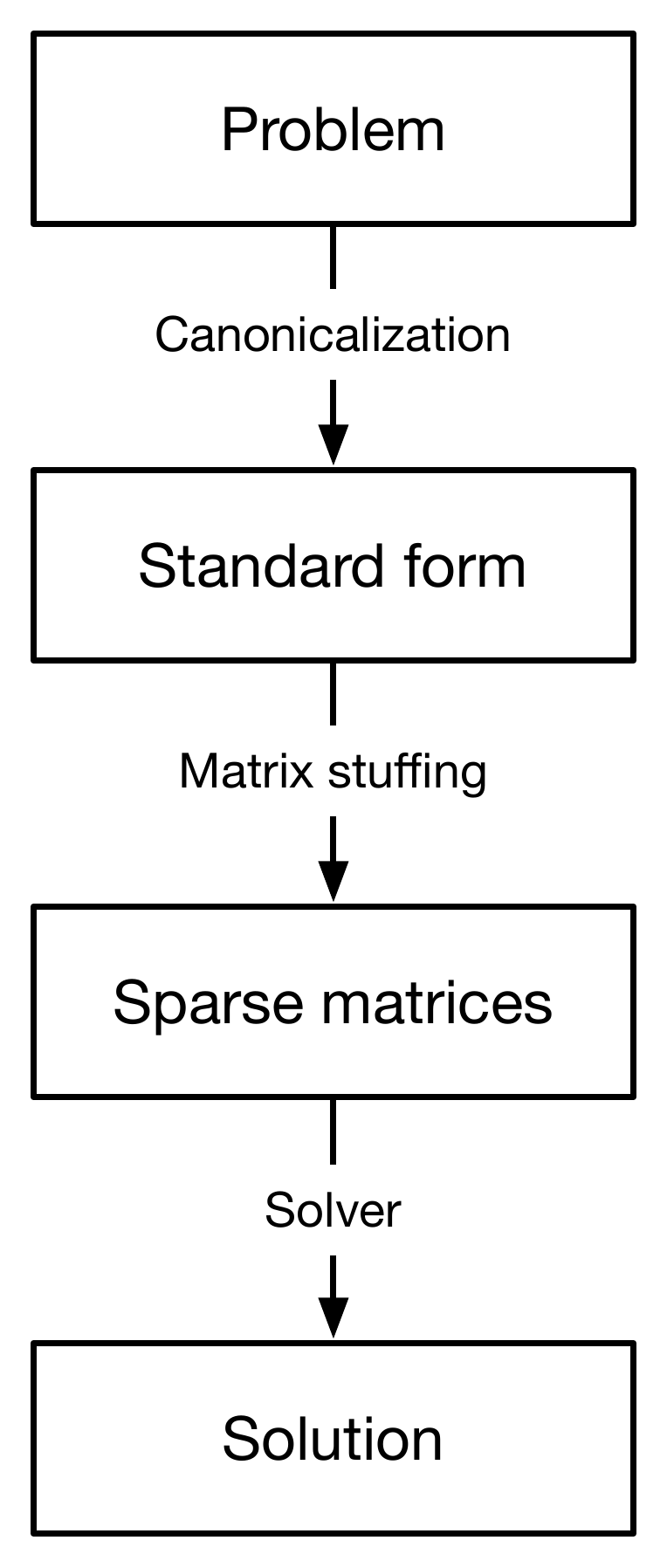}
  \caption{The traditional architecture for optimization modeling frameworks.}
\label{fig-old}
\end{figure}

The traditional architecture for optimization modeling frameworks
dates back to AMPL \cite{AMPL} and GAMS \cite{GAMS} in the 1980s.
In this architecture, solving an optimization problem is divided into a three
step process, shown in Fig.~\ref{fig-old}.
The process begins with a high-level description of the optimization
problem expressed in a modeling language.
The first step is canonicalization,
in which the problem is transformed through symbolic manipulation into
an equivalent problem in a standard form.
The second step is matrix stuffing, in which the symbolic representation
of the standard form is instantiated so that linear operators are represented by
sparse matrices.
Often canonicalization and matrix stuffing are combined into a single step.
The final step is to call a solver with the sparse matrix representation of
the standard form as input
and return the solver output as the solution.

The ecosystem of modeling frameworks for convex optimization is
an illustrative example of the traditional architecture.
Convex optimization modeling languages are built around the
principles of disciplined convex programming,
a set of rules for constructing optimization problems
that make it easy to verify problem convexity.
Implementations include CVX \cite{cvx} and YALMIP \cite{YALMIP} in MATLAB,
CVXPY \cite{diamond2016cvxpy} in Python,
Convex.jl \cite{cvxjl} in Julia,
and the standalone compilers CVXGEN \cite{MB:12} and QCML \cite{QCML}.
We discuss each component of the traditional architecture in the concrete
case of convex modeling frameworks.

\subsection{Canonicalization}
The standard form for convex optimization problems is a
cone program, an optimization problem of the form
\begin{equation}\label{prob-cone}
\begin{array}{ll}
  \mbox{minimize} & c^Tx \\
  \mbox{subject to} & Ax + b \in \mathcal{K},
\end{array}
\end{equation}
where $x \in \reals^n$ is the optimization variable; $A \in \reals^{m \times
  n}$, $b \in \reals^m$, and $c \in \reals^n$ are constants; and
$\mathcal{K}$ is a nonempty closed convex cone \cite{NesNem:92}.
Convex optimization modeling frameworks symbolically convert
problems into cone programs via epigraph transformations \cite{GB:08}.

\subsection{Matrix stuffing}
In solvers and other software that use the cone program standard
form as an input format for problems,
$c$ and $b$ are represented by arrays and $A$ is represented
by a standard sparse matrix format, such as column compressed storage.
Matrix stuffing generates a sparse matrix representation of $A$ from the
symbolic representation generated through canonicalization \cite{DB:15}.
Solvers use a sparse matrix representation of $A$ because they
generally use algorithms and libraries that exploit sparsity.

\subsection{Solver}
A wide variety of solver implementations have been developed for
problems in the cone program standard form.
Many solvers are written in pure C, including MOSEK \cite{mosek}, SDPA \cite{SDPA}, ECOS
\cite{bib:Domahidi2013ecos}, and SCS \cite{SCSpaper}.
Other solvers are written in higher level languages, such as
SeDuMi \cite{Sturm1999} and SDPT3 \cite{toh1999sdpt3} in MATLAB
and CVXOPT \cite{CVXOPT} in Python.
The solvers rely heavily on low-level linear
algebra interfaces like BLAS and LAPACK \cite{lawson1979basic} for basic
operations
and libraries like SuiteSparse \cite{suitesparse} for sparse matrix factorization.
Existing cone solvers are almost exclusively restricted to CPU implementations; an exception is SCS which provides GPU support using the
cuBLAS library \cite{nvidia2008cublas}.

\subsection{Drawbacks}
The traditional approach to optimization modeling frameworks has been
enormously successful, allowing modeling languages and solver implementations to
be developed independently in the programming languages best suited to their
function. The conventional solver implementation is based on interior point
methods, for which the dominant computational effort is solving a sparse linear
system. Such a solver can be ported relatively easily to new platforms provided
the necessary linear algebra libraries (BLAS, LAPACK, SuiteSparse, \etc) are available.

However, many optimization problems of interest are too large to be solved with
interior point methods and, more generally, any method that requires a direct solution
to a linear system involving the $A$ matrix of the cone program standard form
\eqref{prob-cone}. In some problem domains the memory requirements even for
sparse $A$ matrices can be prohibitive (\eg, 2D convolution in large-scale image
reconstruction), while at the same time efficient procedural evaluations of the
matrix-vector computations with $A$ and $A^T$ exist (\eg, FFT-based
convolution).

A possible solution is a first-order method, such as SCS \cite{SCSpaper}, which
only requires solving linear systems to moderate accuracy. This approach can be
implemented with either a direct or indirect method for the linear solver
subroutine. In the case of a direct solver, the computational cost can be
amortized by caching the factorization of the $A$ matrix leading to iterations
that are significantly faster than interior point methods. In the indirect case,
a matrix-free method such as conjugate gradient is used, requiring only
matrix-vector computations with $A$ and $A^T$. In the traditional architecture,
these computations are simply implemented with sparse matrix multiplies, but the
proposed graph-based approach enables taking advantage of specialized linear
operator implementations, as we will discuss in detail in the next section.

\section{Alternative approaches}\label{sec-alternative}


Prior work has explored alternative approaches to bypassing the limitations of
the traditional architecture for optimization modeling frameworks, with a
focus on scaling to larger problem sizes. There are two main lines of work that
are precursors to the graph-based architecture proposed in this paper: the first
replaces the sparse matrix representation of the standard form generated by
matrix stuffing with a more general representation,
and the second explores new standard forms based on functions with
efficient proximal operators. In this section, we review these approaches,
providing motivation for our general graph-based framework.

\subsection{Abstract linear operators}
In solving many convex optimization problems, the majority of computational time is spent in
evaluating linear operators. While the sparse matrix representation of cone
programs is fully general,
it does not provide the most efficient implementation for many types of linear
functions. Matrix-free CVXPY replaces traditional sparse matrix representation of the cone program standard form
\eqref{prob-cone} with a computation graph based representation.
The computation graph representation allows the modeling layer to encode information
about structured linear operators in the optimization problem that solvers can
exploit \cite{DB:15}. The matrix-free CVXPY implementation includes a custom runtime
system for computation graphs, as opposed to the \texttt{cvxflow} implementation
presented in this paper, which is built on TensorFlow.

\subsection{Proximal standard forms}
Another line of work explores solvers based on functions with efficient proximal
operators. Epsilon \cite{wytock2015convex} introduces the standard
form
\begin{equation}\label{prox-prob}
\begin{array}{ll}
  \mbox{minimize} & \sum_{i=1}^N f_i(A_i x),
\end{array}
\end{equation}
where $x \in \reals^n$ is the optimization variable, $A_i \in \reals^{m_i \times
  n}$ are linear operators,
and $f_i$ are functions with efficient proximal operators
\cite{parikh2013proximal}.
Epsilon exploits the flexibility of the standard form \eqref{prox-prob}
to rewrite the problem so it can be solved efficiently by a variant of
the alternating direction method of multipliers (ADMM) \cite{BP:11}. Along
similar lines, POGS \cite{fougner2015parameter} introduces a slightly different
standard form, again based on functions with efficient proximal operators, and
includes a highly efficient GPU implementation of an ADMM-based algorithm.

The ProxImaL modeling framework also targets the standard form
\eqref{prox-prob}, but supports a variety of solver algorithms and applies
problem rewritings specialized to optimization problems in imaging
\cite{proximal2016}. ProxImaL moves towards platform independence by generating
solver implementations using Halide \cite{halide2013}. Halide is a language and
compiler that allows for platform independent abstraction of individual
mathematical operations, but not of full algorithms composed of many operations
inside control logic.

These new proximal standard forms are not necessarily incompatible with the
traditional architecture based on sparse matrices. However, as opposed to cone
solvers and in particular interior point methods, the implementation of
algorithms operating on the proximal standard forms is less reliant on sparse
linear algebra and thus there is less benefit from building on existing
sparse linear algebra libraries. These approaches instead require a library of
proximal operator implementations which can benefit greatly from being built on
a high-level framework such as Halide or TensorFlow, providing platform
independence and a highly optimized runtime system.

\section{Graph-based architecture}\label{sec-proposed}
\begin{figure}
  \centering
  \includegraphics[scale=0.5]{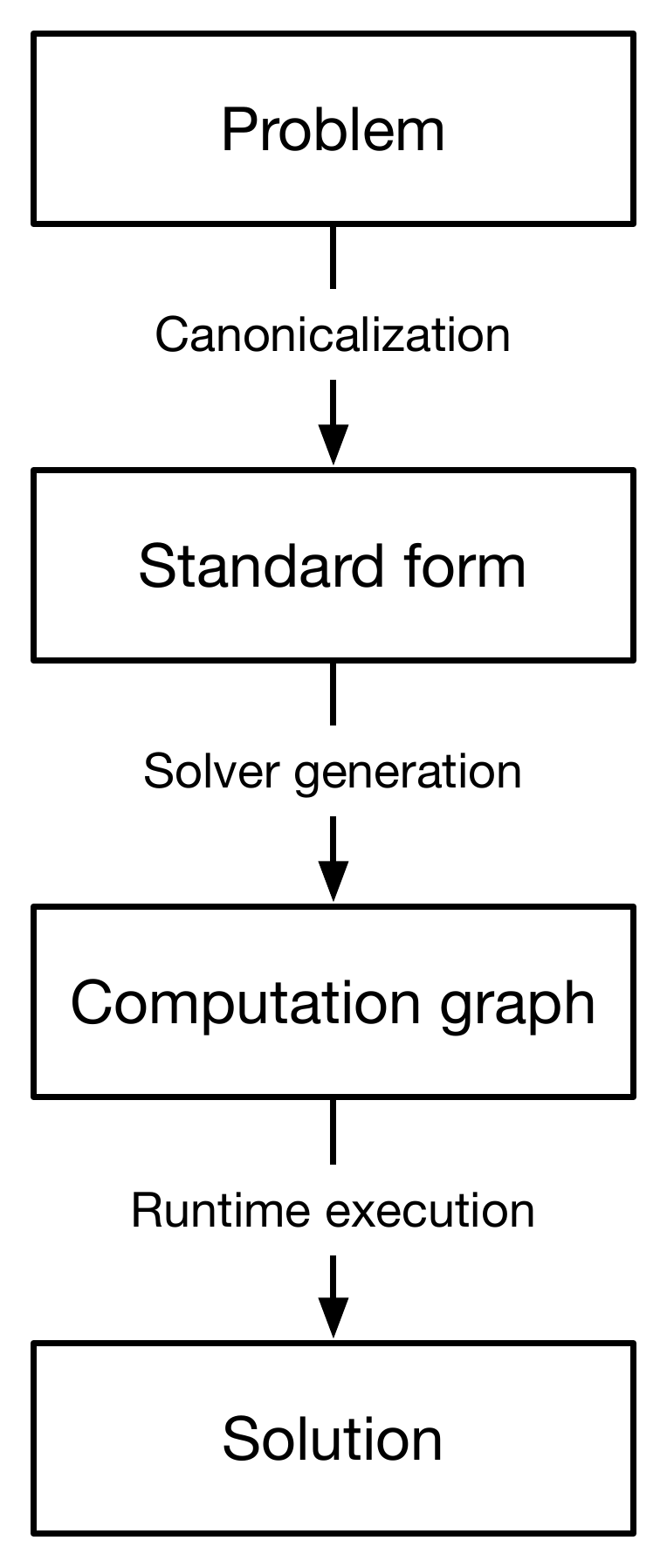}
  \caption{The proposed graph-based architecture for optimization modeling
    frameworks.}
  \label{fig-new}
\end{figure}

In this section, we propose a new graph-based architecture for
optimization modeling frameworks. Our architecture divides the process of
solving an optimization problem into three steps, shown in
Fig.~\ref{fig-new}. As with the traditional architecture we begin with a
high-level problem description which is canonicalized to a standard form. The
solver generation step produces a computation graph representing the solver
algorithm, which is executed by the runtime system to produce a solution.

The key difference from the traditional architecture is that the graph-based
approach directly generates computation graphs representing the numerical
algorithms for solving problems rather than representing all problems with
sparse matrices. The first benefit of this approach is support for abstract linear
operators with highly efficient implementations, such as convolution,
Kronecker products, and others. The second benefit is a closer connection between
canonicalization and solver generation, which can now both be implemented in the
same high-level language and even in a single library. This more easily allows
for supporting different standard forms that incorporate problem-specific
structure. Finally, the new architecture severs the link between the solver
implementation and computing platform, allowing solvers to take advantage of new
computing platforms simply by changing the target of the computation graph runtime system.

We next explain the central abstraction, computation graphs, and describe
how such graphs representing solvers are generated.

\subsection{Computation graphs}

A computation graph is a directed acyclic graph (DAG)
where each vertex represents a mathematical operation
and each edge represents data transfer.
Input vertices have no incoming edges,
while output vertices have no outgoing edges.
A vertex is evaluated by applying its operation to the data
on the vertex and broadcasting the result on its outgoing edges.
The overall graph is evaluated by loading data onto the input vertices,
evaluating the vertices in topological order,
and reading the results off the output vertices.

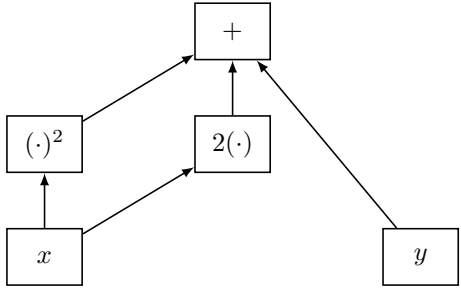
\begin{figure}
\centering
\begin {tikzpicture}[-latex ,auto ,node distance =1.5 cm and 2.5cm ,on grid ,
semithick ,
state/.style ={ rectangle ,top color =white , bottom color = white,
draw,black , text=black , minimum width =1 cm, minimum height = 0.75cm}]
\node[state] (C) {$x$};
\node[state] (A) [above =of C] {$(\cdot)^2$};
\node[state] (B) [above right =of C] {$2(\cdot)$};
\node[state] (E) [below right=of B] {$y$};
\node[state] (D) [above right=of A] {$+$};
\path (C) edge node[below =0.15 cm] {} (A);
\path (C) edge node[below =0.15 cm] {} (B);
\path (A) edge  node[above] {} (D);
\path (B) edge  node[below =0.15 cm] {} (D);
\path (E) edge  node[below =0.15 cm] {} (D);
\end{tikzpicture}
\caption{A computation graph for $f(x,y) = x^2 + 2x + y$.}\label{fig-dag}
\end{figure}

For example, Fig.~\ref{fig-dag} shows a computation graph for
the function $f(x,y) = x^2 + 2x + y$.
The input vertices represent the variables $x$ and $y$.
The output vertex represent the top level sum.
The internal vertices represent the operations $z \to z^2$
and $z \to 2z$.

Given a computation graph to evaluate a function,
computation graphs for evaluating the function's gradient
or adjoint (for linear functions) can be obtained through
simple graph transformations \cite{G:89,DB:15}.
Function, gradient, and adjoint evaluations are the key operations
in first-order and indirect solvers
and are even sufficient to precondition a problem \cite{mf_equil}.

Computation graphs are a useful intermediate representation for solvers
because they abstract away the platform-specific details of both computation
and memory management.
These details are handled by a computation graph runtime system,
which has platform-specific code to execute each mathematical operation
and to pass data from one operation to the next.
By contrast, a solver built on the traditional abstraction of a
low-level linear algebra interface must implement its own platform-specific
logic for mathematical operations not expressible as linear functions and for memory management.

\subsection{Solver generation}

The solver generation step produces a computation graph representing a
numerical algorithm for solving an optimization problem. Graph generation is
naturally implemented in a high-level functional programming style with modular
functions that produce computation graphs implementing numerical
algorithms or subroutines. Typically, these functions take as inputs individual
nodes or in some cases are naturally parameterized by graph generator functions.

As a concrete example, the Python code snippet for generating a TensorFlow graph
representing  the conjugate gradient method for the linear system $Ax = b$ is
shown below.
\begin{verbatim}
def cg_solve(A, b, x_init, tol=1e-8):
  delta = tol*norm(b)

  def body(x, k, r_norm_sq, r, p):
    Ap = A(p)
    alpha = r_norm_sq / dot(p, Ap)
    x = x + alpha*p
    r = r - alpha*Ap
    r_norm_sq_prev = r_norm_sq
    r_norm_sq = dot(r,r)
    beta = r_norm_sq / r_norm_sq_prev
    p = r + beta*p
    return (x, k+1, r_norm_sq, r, p)

  def cond(x, k, r_norm_sq, r, p):
    return tf.sqrt(r_norm_sq) > delta

  r = b - A(x_init)
  loop_vars = (
      x_init, tf.constant(0),
      dot(r, r), r, r)
  return tf.while_loop(
      cond, body, loop_vars)[:3]
\end{verbatim}
In this example, the function \texttt{cg\_solve} is parameterized by the
the linear operator \texttt{A}, and vector \texttt{b} with initial starting
point, vector \texttt{x\_init}. The inputs \texttt{b} and \texttt{x\_init} are
computation graph nodes and \texttt{A} is a single-argument function such that
\texttt{A(x)} produces the computation graph representing the linear operator
applied to an arbitrary vector \texttt{x}. Implemented in this fashion, the
conjugate gradient method can be applied to any linear operator expressed as a
computation graph.

\section{Numerical examples}\label{sec-examples}

In this section, we present numerical examples of solving convex optimization
problems in our proposed architecture. As solving linear systems forms the basis for
convex methods, we first present results for an indirect linear solver with various
linear operators.
Using this indirect linear solver as a subroutine, we then
implement a version of SCS \cite{SCSpaper} in the computation graph framework
and compare with the native version of SCS implemented in C. We present results
for both CPU and GPU environments; all experiments are run on a 32-core Intel
Xeon 2.2GHz processor and an nVidia Titan X GPU with 12GB of RAM.

Our implementation builds on CVXPY \cite{diamond2016cvxpy},
a convex optimization modeling framework in Python.
Using this framework, convex optimization problems can be expressed with minimal code
and are automatically converted into the standard conic form \eqref{prob-cone}.
As an example, the nonnegative deconvolution problem we
consider in Section \ref{sec-nn-deconv} is written as the following Python code.
\begin{verbatim}
from cvxpy import *
x = Variable(n)
f = norm(conv(c, x) - b, 2)
prob = Problem(Minimize(f), [x >= 0])
\end{verbatim}
Here \texttt{c} and \texttt{b} are previously-defined problem inputs
and \texttt{n} is the size of the optimization variable. Our
implementation differs from the existing CVXPY functionality in that instead of
solving problems by constructing sparse matrices and calling numerical routines
written in C, we build a computation graph, as described in Section
\ref{sec-proposed}, and evaluate with TensorFlow. Ultimately, this
implementation achieves faster running times than existing methods---for
example, on the large nonnegative deconvolution example, our implementation takes
roughly 1/10th the time of SCS running on GPU, the
existing state-of-the-art method for solving large convex problems to moderate
accuracy.

Concurrent with the publication of this paper, we are releasing the \texttt{cvxflow}
Python library; it is available at
\begin{quote}
  \url{http://github.com/cvxgrp/cvxflow}
\end{quote}
and includes the code for all of the examples in this section. The implementation
is general and solves any problem modeled with CVXPY using TensorFlow.

\subsection{Regularized least squares}
\label{sec-linear}

We begin with solving linear systems using the conjugate gradient method (CG)
\cite{HeS:52}. CG is matrix-free which makes it a natural fit for
linear systems represented as a graph, allowing for specialized implementations
of each linear operator including those that are inefficient to represent as
sparse matrices such as convolution, Kronecker products, and others. In terms of
the graph-based architecture shown in Fig.~\ref{fig-new}, the standard form in
this example is a linear system and the solver generation step generates a graph
representing the conjugate gradient method.


In particular, we consider the regularized least squares problem
\begin{equation}
  \begin{array}{ll}
    \mbox{minimize} & (1/2)\|Ax - b\|_2^2 + \lambda \|x\|_2^2
  \end{array}
\end{equation}
where $x \in \reals^n$ is the optimization variable, the linear map $A :
\reals^n \to \reals^m$ and vector $b \in \reals^m$ are problem data, and $\lambda
> 0$ is the regularization parameter. This problem has the solution
\begin{equation}
  x^\star = (\lambda I + A^TA)^{-1}A^Tb,
\end{equation}
which can be computed by solving a linear system.

It is often the case that $A$ takes the form of a sparse or
dense matrix; for example, in a statistical problem each row of $A$
may represent an observation of multiple variables weighted by $x$ in order to
predict the response variable. However, $A$ can also
be an abstract linear operator; for example, a convolution with a vector $c$,
written as $Ax = c * x$. We present results for each of these examples: a sparse matrix,
a dense matrix, and convolution.

In the matrix examples, entries are sampled
from $\mathcal{N}(0,1)$ with 1\% nonzero in the sparse case. For convolution, we
apply the Gaussian kernel with standard deviation $n/10$. In all cases, the
response variable is formed by $ b = A \hat{x} + v$ where $v$
has entries sampled from $\mathcal{N}(0, 0.01^2)$ and $\hat{x}$ from
$\mathcal{N}(0,1)$. The conjugate
gradient method is run until the residual satisfies $\|(\lambda I + A^TA)x^k -
A^Tb\|_2 \le  10^{-8} \|A^Tb\|_2$.

\begin{table}
  \centering
  \begin{tabular}{llll}
                            & \textbf{dense matrix} & \textbf{sparse matrix} & \textbf{convolution} \\
    \hline
    variables $n$           & 3000                   & 3000                    & 3000 \\
    nonzeros in $A$         & 18000000               & 180000                  & 4095000 \\
    \\ \textbf{spsolve} \\
    solve time              & 255 secs              & 28 secs               & 41 secs \\
    memory usage            & 2.2 GB                & 1.06 GB                 & 1.5 GB \\
    objective               & $5.97 \times 10^{-1}$ & $5.97 \times 10^{-1}$ & $7.68 \times 10^{-1}$ \\
    \\
    \textbf{CG TensorFlow} \\
    solve time, CPU         & 3.0 secs             & 0.9 secs              & 2.9 secs \\
    solve time, GPU         & 2.0 secs             & 0.7 secs              & 1.0 secs\\
    graph build time        & 0.4 secs             & 0.1 secs              & 0.1 secs \\
    memory usage            & 1.8 GB               & 755 MB                 & 946 MB \\
    objective               & $5.97 \times 10^{1}$ & $5.97 \times 10^{-1}$ & $7.68 \times 10^{-1}$ \\
    CG iterations           & 49                   & 49                     & 71 \\
  \end{tabular}
  \caption{Results for regularized least squares.}
  \label{tab-ls}
\end{table}

Table \ref{tab-ls} shows the results for these experiments, demonstrating that
conjugate gradient on TensorFlow is significantly faster than the baseline method,
\texttt{scipy.sparse.spsolve}. This is a somewhat weak baseline as
\texttt{spsolve} does not run on GPU  and is not
well-suited for dense matrices. Nonetheless, this comparison highlights the
difference in architecture exploited by TensorFlow which can take advantage of
dedicated implementations for the linear operators leading to significantly
faster solve times.


\subsection{Lasso}

\begin{table}
  \centering
  \begin{tabular}{llll}
                       & \textbf{dense matrix}      & \textbf{sparse matrix}     & \textbf{convolution} \\
    \hline
    variables $n$      & 6001                 & 6001                 & 6001      \\
    constraints $m$    & 12002                & 12002                & 12001     \\
    nonzeros in $A$    & 18012002             & 1812002              & 4107002   \\
    \\ \textbf{SCS native} \\
    solve time, CPU    & 29 secs           & 3.4 secs            & 6.4 secs \\
    solve time, GPU    & 27 secs           & 3.8 secs            & 7.6 secs \\
    matrix build time  & 13 secs           & 1.4 secs            & 2.8 secs \\
    memory usage       & 3.1 GB              & 663 MB              & 927 MB \\
    objective          & $3.36 \times 10^1$  & $3.19 \times 10^1$   & $2.02 \times 10^0$ \\
    SCS iterations         & 40                  & 40                 & 60            \\
    avg. CG iterations & 2.66                & 2.71               & 2.72           \\
    \\ \textbf{SCS TensorFlow} \\
    solve time, CPU    & 23 secs           & 25 secs          & 24 secs \\
    solve time, GPU    & 9.9 secs            & 7.1 secs           & 5.3 secs \\
    graph build time   & 1.8 secs            & 2.0 secs           & 0.8 secs    \\
    memory usage       & 8.7 GB              & 4.6 GB             & 1.2 GB \\
    objective          & $3.36 \times 10^1$  & $3.19 \times 10^1$ & $2.02 \times 10^0$ \\
    SCS iterations         & 60                  & 40                 & 180 \\
    avg. CG iterations & 3.35                & 3.55               & 1.93 \\
  \end{tabular}
  \caption{Results for lasso.}
  \label{tab-lasso}
\end{table}

Next we solve a convex problem with SCS \cite{SCSpaper}. In this
case, the canonicalization step produces a problem in the standard cone form
\eqref{prob-cone} and solver generation produces a graph implementing the SCS
iterations. In essence, the algorithm iterates
between projections onto a linear subspace and a convex cone; the
former is done through solving a linear system with a computation graph
representing the CG method as in the previous section.
The SCS method is appealing in this context as it works well with approximate
solutions to linear systems, such as those produced by CG.

We consider the lasso problem
\begin{equation}
  \begin{array}{ll}
    \mbox{minimize} & (1/2)\|Ax - b\|_2^2 + \lambda \|x\|_1,
  \end{array}
\end{equation}
where the regularization term $\|x\|_1$ replaces the $\|x\|_2^2$ in the
regularized least squares problem from the previous section. This problem is
convex but no longer has a closed-form solution.

To generate problem instances, we construct example linear operators $A$ as in the
previous section. We set the regularization parameter to
$\lambda = 0.1 \|A^Tb\|_\infty$ where $\|A^Tb\|_\infty$
is the smallest value of $\lambda$ such that the solution is zero.

Table \ref{tab-lasso} compares the TensorFlow version of SCS to the native implementation
and demonstrates that in the dense matrix and convolution cases, the solve time
on GPU is faster with TensorFlow. This highlights the benefit of
the computation graph, taking advantage of specialized implementations for dense
matrix multiplication and convolution. In contrast, when the input linear
operator $A$ is a sparse matrix, native SCS is faster.

\subsection{Nonnegative deconvolution}
\label{sec-nn-deconv}

As a final example further illustrating the benefit of abstract linear
operators, we consider the nonnegative deconvolution problem
\begin{equation}
  \begin{array}{ll}
    \mbox{minimize} & \|c * x - b \|_2 \\
    \mbox{subject to} & x \ge 0
  \end{array}
\end{equation}
where $x \in \reals^n$ is the optimization variable, and $c \in \reals^n$, $b \in
\reals^{2n-1}$ are problem data. As in the previous example, the
canonicalization step transforms the problem to the
standard form \eqref{prob-cone} and solver generation produces a computation
graph for the SCS algorithm.

We generate problem instances by taking $c$ to be
the Gaussian kernel with standard deviation $n/10$ and convolving it with a
sparse signal $\hat{x}$ with 5 nonzero entries sampled
uniformly from $[0,n/10]$. We set the response $b = c * \hat{x} + v$ with $v \sim
\mathcal{N}(0, 0.01^2)$.

\begin{table}
  \centering
  \begin{tabular}{llll}
                       & \textbf{small}      & \textbf{medium}     & \textbf{large} \\
    \hline
    variables $n$      & 101                 & 1001                & 10001     \\
    constraints $m$    & 300                 & 3000                & 30000       \\
    nonzeros in $A$    & 9401                & 816001              & 69220001   \\
    \\ \textbf{SCS native} \\
    solve time, CPU    & 0.1 secs           & 2.2 secs           & 260 secs \\
    solve time, GPU    & 2.0 secs           & 2.0 secs            & 105 secs \\
    matrix build time  & 0.01 secs           & 0.6 secs            & 52 secs \\
    memory usage       & 360 MB              & 470 MB              & 10.4 GB \\
    objective          & $1.38 \times 10^0$ & $4.57 \times 10^0$   & $1.41 \times 10^1$ \\
    SCS iterations         & 380                 & 100                 & 160            \\
    avg. CG iterations & 8.44                & 2.95               & 3.01           \\
    \\ \textbf{SCS TensorFlow} \\
    solve time, CPU    & 3.4 secs           & 5.7 secs          & 88 secs \\
    solve time, GPU    & 5.7 secs           & 3.2 secs           & 13 secs \\
    graph build time   & 0.8 secs           & 0.8 secs            & 0.9 secs    \\
    memory usage       & 895 MB              & 984 MB             & 1.3 GB \\
    objective          & $1.38 \times 10^0$  & $4.57 \times 10^0$  & $1.41 \times 10^1$ \\
    SCS iterations         & 480                 & 100                 & 160 \\
    avg. CG iterations & 2.75                & 2.00                & 2.00 \\
  \end{tabular}
  \caption{Results for nonnegative deconvolution.}
  \label{tab-nn-deconv}
\end{table}

Table \ref{tab-nn-deconv} shows that on large problem
sizes, the SCS TensorFlow implementation performs significantly better
than the native implementation, requiring 13 seconds as compared to 105
seconds. This difference is largely due to differences
in architecture, as the matrix-based SCS requires a considerable amount of time
(52 seconds) to simply construct the sparse matrix representing the convolution
operator. As many linear operators benefit from from specialized implementations
(see \eg, \cite{SPOT,BEFB,VaB:95,diamond2015convex}),
one could easily demonstrate an even more significant
gap between the proposed architecture and existing methods simply by choosing
more egregious examples that highlight this difference.

\section*{Acknowledgments}
This material is based upon work supported by the
National Science Foundation Graduate
Research Fellowship under Grant No. DGE-114747
and by the DARPA XDATA program.

\newpage

\vskip 0.2in
\bibliography{cvxflow}

\end{document}